\documentclass[12pt]{amsart}
\usepackage{amssymb}

\usepackage{t1enc}
\usepackage[latin2]{inputenc}
\usepackage{verbatim}
\usepackage{amsmath,amsfonts,amssymb,amsthm}
\usepackage[mathcal]{eucal}
\usepackage{enumerate}
\usepackage[centertags]{amsmath}
\usepackage{graphics}

\newtheorem{theorem}{Theorem}
\newtheorem{lemma}{Lemma}

\newtheorem{corollary}{Corollary}

\numberwithin{equation}{subsection}

\begin{document}
\author{George Tephnadze}
\title[Fej\'er means ]{Fej\'er means of Vilenkin-Fourier series}
\address{G. Tephnadze, Department of Mathematics, Faculty of Exact and Natural
Sciences, Tbilisi State University, Chavchavadze str. 1, Tbilisi 0128,
Georgia}
\email{giorgitephnadze@gmail.com}
\date{}
\maketitle

\begin{abstract}
The main aim of this paper is to prove that there exist a martingale $f\in
H_{1/2}$ such that Fej\'er means of Vilenkin-Fourier series of the
martingale $f$ is not uniformly bounded in the space $L_{1/2}.$
\end{abstract}

\textbf{2000 Mathematics Subject Classification.} 42C10.

\textbf{Key words and phrases:} Vilenkin system, Fejér means, martingale
Hardy space.

\section{ INTRODUCTION}

In one-dimensional case the weak type inequality

\begin{equation*}
\mu \left( \sigma ^{*}f>\lambda \right) \leq \frac{c}{\lambda }\left\|
f\right\| _{1},\text{ \qquad }\left( \lambda >0\right)
\end{equation*}
can be found in Zygmund \cite{Zy} for the trigonometric series, in Schipp
\cite{Sc} for Walsh series and in Pál, Simon \cite{PS} for bounded Vilenkin
series. Again in one-dimensional, Fujji \cite{Fu} and Simon \cite{Si2}
verified that $\sigma ^{*}$ is bounded from $H_{1}$ to $L_{1}$. Weisz \cite
{We2} generalized this result and proved the boundedness of $\sigma ^{*}$
from the martingale space $H_{p}$ to the space $L_{p}$ for $p>1/2$. Simon
\cite{Si1} gave a counterexample, which shows that boundedness does not hold
for $0<p<1/2.$ The counterexample for $p=1/2$ due to Goginava ( \cite{Go},
see also \cite{BGG2}).

In \cite{BGG2} the following is proved:

For any bounded Vilenkin system the maximal operator of the Fejér means is
not bounded from the martingale Hardy space $H_{1/2}$ to the space $L_{1/2}.$

In this paper we shall prove a stronger result then the unboundednes of the
maximal operator from the Hardy space $H_{1/2}$ to the space $L_{1/2},$ in
particular, we shall prove that there exists a martingale $f\in H_{1/2}$
such that Fejér means of Vilenkin-Fourier series of the martingale $f$ is
not uniformly bounded in the the space $L_{1/2}.$

\section{DEFINITONS AND NOTATIONS}

Let $N_{+}$ denote the set of the positive integers, $N:=N_{+}\cup \{0\}.$
Let $m:=(m_{0,}m_{1....})$ denote a sequence of the positive integers, not
less than 2. Denote by $Z_{m_{k}}:=\{0,1,...m_{k}-1\}$ the addition group of
integers modulo $m_{k}$.

Define the group $G_{m}$ as the complete direct product of the groups $%
Z_{m_{i}},$ with the product of the discrete topologies of $Z_{m_{j}}$ $^{,}$%
s.

The direct product $\mu $ of the measures

\begin{equation*}
\mu _{k}\left( \{j\}\right) :=1/m_{k}\text{,\qquad }(j\in Z_{m_{k}})
\end{equation*}
is the Haar measure on $G_{m_{k\text{ }}}$with $\mu \left( G_{m}\right) =1.$

If $\sup\limits_{n}m_{n}<\infty $, then we call $G_{m}$ a bounded Vilenkin
group. If the generating sequence $m$ is not bounded, then $G_{m}$ is said
to be an unbounded Vilenkin group. \textbf{In this paper we discuss bounded
Vilenkin groups only.}

The elements of $G_{m}$ represented by sequences

\begin{equation*}
x:=\left( x_{0},x_{1},...,x_{j},...\right) ,\text{ \qquad }\left( x_{i}\in
Z_{m_{j}}\right) .
\end{equation*}

It is easy to give a base for the neighborhood of $G_{m}:$

\begin{equation*}
I_{0}\left( x\right) :=G_{m,}
\end{equation*}

\begin{equation*}
I_{n}(x):=\{y\in G_{m}\mid y_{0}=x_{0},...y_{n-1}=x_{n-1}\},\,\,\left( x\in
G_{m},\text{ }n\in N\right) .
\end{equation*}
Denote $I_{n}:=I_{n}\left( 0\right) ,$ for $n\in N_{+}.$

If we define the so-called generalized number system based on $m$ in the
following way :

\begin{equation*}
M_{0}:=1,\text{ \qquad }M_{k+1}:=m_{k}M_{k},\,\,\,(k\in N),
\end{equation*}
then every $n\in N$ can be uniquely expressed as $n=\sum_{j=0}^{\infty
}n_{j}M_{j},$ where $n_{j}\in Z_{m_{j}},$ $(j\in N_{+})$ and only a finite
number of $n_{j}`$s differ from zero.

Next, we introduce on $G_{m}$ an ortonormal system, which is called the
Vilenkin system. At first define the complex valued function $r_{k}\left(
x\right) :G_{m}\rightarrow C,$ The generalized Rademacher functions as

\begin{equation*}
r_{k}\left( x\right) :=\exp \left( 2\pi ix_{k}/m_{k}\right) ,\text{ \qquad }%
\left( i^{2}=-1,\text{ }x\in G_{m},\text{ }k\in N\right) .
\end{equation*}

Now define the Vilenkin system$\,\,\,\psi :=(\psi _{n}:n\in N)$ on $G_{m}$
as:

\begin{equation*}
\psi _{n}(x):=\prod\limits_{k=0}^{\infty }r_{k}^{n_{k}}\left( x\right)
,\,\,\,\left( n\in N\right) .
\end{equation*}
Specifically, we call this system the Walsh-Paley one if $m\equiv 2.$

The Vilenkin system is orthonormal and complete in $L_{2}\left( G_{m}\right)
$ \cite{AVD,Vi}.

Now we introduce analogues of the usual definitions in Fourier-analysis. If $%
f\in L_{1}\left( G_{m}\right) $ we can establish the Fourier coefficients,
the partial sums of the Fourier series, the Fejér means, the Dirichlet
kernels with respect to the Vilenkin system $\psi $ in the usual manner:

$\qquad $%
\begin{eqnarray*}
\widehat{f}\left( k\right) &:&=\int_{G_{m}}f\overline{\psi }_{k}d\mu ,\text{
\qquad }\left( k\in N\right) , \\
S_{n}f &:&=\sum_{k=0}^{n-1}\widehat{f}\left( k\right) \psi _{k},\text{
\qquad }\left( n\in N_{+},\text{ }S_{0}f:=0\right) , \\
\sigma _{n}f &:&=\frac{1}{n}\sum_{k=0}^{n-1}S_{k}f\,\,,\text{ \qquad }%
\,\,\left( n\in N_{+}\right) , \\
D_{n} &:&=\sum_{k=0}^{n-1}\psi _{k\text{ }},\,\qquad \qquad \,\,\left( n\in
N_{+}\right) .
\end{eqnarray*}

Recall that
\begin{equation*}
D_{M_{n}}\left( x\right) =\left\{
\begin{array}{l}
M_{n},\text{ \qquad }\,\,\text{if\thinspace \thinspace \thinspace }x\in
I_{n}, \\
0,\,\qquad \,\text{if}\,\,x\notin I_{n}.
\end{array}
\right.
\end{equation*}

The norm (or quasinorm) of the space $L_{p}(G_{m})$ is defined by \qquad

\begin{equation*}
\left\| f\right\| _{p}:=\left( \int_{G_{m}}\left| f(x)\right| ^{p}d\mu
(x)\right) ^{\frac{1}{p}},\qquad \left( 0<p<\infty \right) .
\end{equation*}

The $\sigma -$algebra generated by the intervals $\left\{ I_{n}\left(
x\right) :x\in G_{m}\right\} $ will be denoted by $\digamma _{n}\left( n\in
N\right) .$ Denote by $f=\left( f^{\left( n\right) },n\in N\right) $ a
martingale with respect to $\digamma _{n}\left( n\in N\right) .$ (for
details see e.g. \cite{We1}).

The maximal function of a martingale $f$ is defined by

\begin{equation*}
f^{*}=\sup_{n\in N}\left| f^{(n)}\right| .
\end{equation*}

In case $f\in L_{1}\left( G_{m}\right) ,$the maximal functions are also be
given by

\begin{equation*}
f^{*}\left( x\right) =\sup\limits_{n\in N}\frac{1}{\mu \left( I_{n}\left(
x\right) \right) }\left| \int\limits_{I_{n}\left( x\right) }f\left( u\right)
d\mu \left( u\right) \right| .
\end{equation*}

For $0<p<\infty ,$ the Hardy martingale spaces $H_{p}$ $\left( G_{m}\right) $
consist of all martingale, for which

\begin{equation*}
\left\| f\right\| _{H_{p}}:=\left\| f^{*}\right\| _{L_{p}}<\infty .
\end{equation*}

If $f\in L_{1}\left( G_{m}\right) ,$ then it is easy to show that the
sequence $\left( S_{M_{n}}\left( f\right) :n\in N\right) $ is a martingale.

If $f=\left( f^{\left( n\right) },n\in N\right) $ is martingale then the
Vilenkin-Fourier coefficients must be defined in a slightly different manner:

\begin{equation*}
\widehat{f}\left( i\right) :=\lim_{k\rightarrow \infty
}\int_{G_{m}}f^{\left( k\right) }\left( x\right) \overline{\Psi }_{i}\left(
x\right) d\mu \left( x\right) .
\end{equation*}

The Vilenkin-Fourier coefficients of $f\in L_{1}\left( G_{m}\right) $ are
the same as those of the martingale $\left( S_{M_{n}}\left( f\right) :n\in
N\right) $ obtained from $f$.

For a martingale $f$ the maximal operators of the Fej\'er means are defined
by $\qquad \qquad $

\begin{equation*}
\sigma ^{*}f\left( x\right) =\sup_{n\in N}\left| \sigma _{n}f(x)\right| .
\end{equation*}

A bounded measurable function $a$ is p-atom, if there exists a interval I,
such that

\begin{equation*}
\left\{
\begin{array}{l}
a)\qquad \int_{I}ad\mu =0, \\
b)\ \qquad \left\| a\right\| _{\infty }\leq \mu \left( I\right) ^{\frac{-1}{p%
}}, \\
c)\qquad \text{supp}\left( a\right) \subset I.\qquad
\end{array}
\right.
\end{equation*}

\section{FORMULATION OF MAIN RESULT}

\begin{theorem}
There exist a martingale $f\in H_{1/2}$ such that
\begin{equation*}
\sup\limits_{n}\left\| \sigma _{n}f\right\| _{1/2}=+\infty .
\end{equation*}
\end{theorem}

\begin{corollary}
There exist a martingale $f\in H_{1/2}$ such that
\begin{equation*}
\left\| \sigma ^{*}f\right\| _{1/2}=+\infty .
\end{equation*}
\end{corollary}

\bigskip

\section{AUXILIARY PROPOSITIONS}

\begin{lemma}
\cite{We3} A martingale $f=\left( f^{\left( n\right) },n\in N\right) $ is in
$H_{p}\left( 0<p\leq 1\right) $ if and only if there exist a sequence $%
\left( a_{k},k\in N\right) $ of p-atoms and a sequence $\left( \mu _{k},k\in
N\right) ,$ of a real numbers, such that for every $n\in N:$
\end{lemma}

\begin{equation}
\qquad \sum_{k=0}^{\infty }\mu _{k}S_{M_{n}}a_{k}=f^{\left( n\right) },
\label{1}
\end{equation}

\begin{equation*}
\qquad \sum_{k=0}^{\infty }\left| \mu _{k}\right| ^{p}<\infty .
\end{equation*}
Moreover, $\left\| f\right\| _{H_{p}}\backsim \inf \left( \sum_{k=0}^{\infty
}\left| \mu _{k}\right| ^{p}\right) ^{1/p}$, where the infimum is taken over
all decomposition of $f$ of the form (\ref{1}).

\begin{lemma}
\cite{BGG} Let $2<A\in N_{+},$ $k\leq s<A$ and $%
q_{A}=M_{2A}+M_{2A-2}+...+M_{2}+M_{0},$ then

\begin{equation*}
q_{A-1}\left| K_{q_{A-1}}(x)\right| \geq \frac{M_{2k}M_{2s}}{4}.
\end{equation*}

for

\begin{equation*}
x\in I_{2A}\left( 0,...,x_{2k}\neq 0,0,...,0,x_{2s}\neq
0,x_{2s+1},...x_{2A-1}\right) ,
\end{equation*}
\end{lemma}

\begin{equation*}
k=0,1,...,A-3,\qquad s=k+2,k+3,...,A-1.
\end{equation*}

\section{PROOF OF THE THEOREM}

Let $\left\{ \alpha _{k}:k\in N\right\} $ be an increasing sequence of the
positive integers such that:

\qquad
\begin{equation}
\sum_{k=0}^{\infty }\alpha _{k}^{-1/2}<\infty ,  \label{2}
\end{equation}

\begin{equation}
\sum_{\eta =0}^{k-1}\frac{\left( M_{2\alpha _{\eta }}\right) ^{2}}{\alpha
_{\eta }}<\frac{\left( M_{2\alpha _{k}}\right) ^{2}}{\alpha _{k}},  \label{3}
\end{equation}

\begin{equation}
\qquad \frac{32M\left( M_{2\alpha _{k-1}}\right) ^{2}}{\alpha _{k-1}}<\frac{%
M_{\alpha _{k}}}{\alpha _{k}},  \label{4}
\end{equation}
where $M=\sup \left\{ m_{0},\text{ }m_{1}\text{ }...\right\} ,$ $\left(
2\leq M<\infty \right) .$

We note that such an increasing sequence $\left\{ \alpha _{k}:k\in N\right\}
$ which satisfies conditions (\ref{2}-\ref{4}) can be constructed.

Let \qquad
\begin{equation*}
f^{\left( A\right) }\left( x\right) =\sum_{\left\{ k;\text{ }2\alpha
_{k}<A\right\} }\lambda _{k}a_{k},
\end{equation*}
where
\begin{equation*}
\lambda _{k}=\frac{1}{\alpha _{k}},
\end{equation*}
and

\begin{equation*}
a_{k}\left( x\right) =\frac{M_{2\alpha _{k}}}{M}\left( D_{M_{(2\alpha
_{k}+1)}}\left( x\right) -D_{M_{_{2\alpha _{k}}}}\left( x\right) \right) .
\end{equation*}

It is easy to show that the martingale $\,f=\left( f^{\left( 1\right)
},f^{\left( 2\right) }...f^{\left( A\right) }...\right) \in H_{1/2}.$

Indeed, since

\begin{equation}
S_{M_{A}}a_{k}\left( x\right) =\left\{
\begin{array}{l}
a_{k}\left( x\right) \text{, \quad }2\alpha _{k}<A, \\
0\text{ , \qquad }2\alpha _{k}\geq A,
\end{array}
\right.  \label{5}
\end{equation}

\begin{eqnarray*}
\text{supp}(a_{k}) &=&I_{2\alpha _{k}}, \\
\int_{I_{2\alpha _{k}}}a_{k}d\mu &=&0
\end{eqnarray*}

and

\begin{equation*}
\left\| a_{k}\right\| _{\infty }\leq \frac{M_{2\alpha _{k}}}{M}M_{2\alpha
_{k}+1}\leq (M_{2\alpha _{k}})^{2}=(\text{supp }a_{k})^{-2}.
\end{equation*}
if we apply lemma 1 and (\ref{2}) we conclude that $f\in H_{1/2}.$

It is easy to show that

\begin{equation}
\widehat{f}(j)=\left\{
\begin{array}{l}
\frac{1}{M}\frac{M_{2\alpha _{k}}}{\cdot \alpha _{k}},\,\,\text{ if
\thinspace \thinspace }j\in \left\{ M_{2\alpha _{k}},...,\text{ ~}M_{2\alpha
_{k}+1}-1\right\} ,\text{ }k=0,1,2..., \\
0\text{ },\text{ \thinspace \qquad \thinspace \thinspace if \thinspace
\thinspace \thinspace }j\notin \bigcup\limits_{k=1}^{\infty }\left\{
M_{2\alpha _{k}},...,\text{ ~}M_{2\alpha _{k}+1}-1\right\} .\text{ }
\end{array}
\right.  \label{6}
\end{equation}

We can write

\begin{equation}
\sigma _{q_{\alpha _{k}}}f\left( x\right) =\frac{1}{q_{\alpha _{k}}}%
\sum_{j=0}^{M_{2\alpha _{k}}-1}S_{j}f\left( x\right) +\frac{1}{q_{\alpha
_{k}}}\sum_{j=M_{2\alpha _{k}}}^{q_{\alpha _{k}}-1}S_{j}f\left( x\right)
=I+II.  \label{7}
\end{equation}
Let $M_{2\alpha _{k}}\leq j<q_{\alpha _{k}}.$ Then applying (\ref{6}) we have

\begin{eqnarray}
&&S_{j}f\left( x\right) =\sum_{v=0}^{M_{2\alpha _{k-1}+1}-1}\widehat{f}%
(v)\psi _{v}\left( x\right) +\sum_{v=M_{2\alpha _{k}}}^{j-1}\widehat{f}%
(v)\psi _{v}\left( x\right)  \label{8} \\
&=&\sum_{\eta =0}^{k-1}\sum_{v=M_{2\alpha _{\eta }}}^{M_{2\alpha _{\eta
}+1}-1}\widehat{f}(v)\psi _{v}\left( x\right) +\sum_{v=M_{2\alpha
_{k}}}^{j-1}\widehat{f}(v)\psi _{v}\left( x\right)  \notag \\
&=&\frac{1}{M}\sum_{\eta =0}^{k-1}\sum_{v=M_{2\alpha _{\eta }}}^{M_{2\alpha
_{\eta }+1}-1}\frac{M_{2\alpha _{\eta }}}{\alpha _{\eta }}\psi _{v}\left(
x\right) +\frac{1}{M}\frac{M_{2\alpha _{k}}}{\alpha _{k}}\sum_{v=M_{2\alpha
_{k}}}^{j-1}\psi _{v}\left( x\right)  \notag \\
&=&\frac{1}{M}\sum_{\eta =0}^{k-1}\frac{M_{2\alpha _{\eta }}}{\alpha _{\eta }%
}\left( D_{M_{_{2\alpha _{\eta }+1}}}\left( x\right) -D_{M_{_{2\alpha _{\eta
}}}}\left( x\right) \right)  \notag \\
&&+\frac{1}{M}\frac{M_{2\alpha _{k}}}{\alpha _{k}}\left( D_{_{j}}\left(
x\right) -D_{M_{_{2\alpha _{k}}}}\left( x\right) \right) .  \notag
\end{eqnarray}

Applying (\ref{8}) in II we have

\begin{eqnarray*}
II &=&\frac{1}{M}\frac{q_{\alpha _{k}}-M_{2\alpha _{k}}}{q_{\alpha _{k}}}%
\sum_{\eta =0}^{k-1}\frac{M_{2\alpha _{\eta }}}{\alpha _{\eta }}\left(
D_{M_{_{2\alpha _{\eta }+1}}}\left( x\right) -D_{M_{_{2\alpha _{\eta
}}}}\left( x\right) \right) \\
&&+\frac{1}{M}\frac{M_{2\alpha _{k}}}{\alpha _{k}q_{\alpha _{k}}}%
\sum_{j=M_{2\alpha _{k}}}^{q_{\alpha _{k}}-1}\left( D_{_{j}}\left( x\right)
-D_{M_{_{2\alpha _{k}}}}\left( x\right) \right) \\
&=&II_{1}+II_{2}.
\end{eqnarray*}

It is evident

\begin{equation*}
\left| \frac{q_{\alpha _{k}}-M_{2\alpha _{k}}}{q_{\alpha _{k}}}\right| <1
\end{equation*}
and

\begin{eqnarray*}
&&\left| \left( D_{M_{_{2\alpha _{\eta }+1}}}\left( x\right)
-D_{M_{_{2\alpha _{\eta }}}}\left( x\right) \right) \right| \\
&\leq &M_{_{2\alpha _{\eta }+1}}=m_{_{2\alpha _{\eta }}}M_{_{2\alpha _{\eta
}}}\leq M\cdot M_{_{2\alpha _{\eta }}}.
\end{eqnarray*}

Applying (\ref{3}) we have

\begin{equation}
\left| II_{1}\right| \leq \sum_{\eta =0}^{k-1}\frac{M_{2\alpha _{\eta }}}{%
\alpha _{\eta }}\frac{1}{M}M\cdot M_{2\alpha _{\eta }}\leq \frac{2\left(
M_{2\alpha _{k-1}}\right) ^{2}}{\alpha _{k-1}}.  \label{9}
\end{equation}

Since

\begin{equation*}
D_{j+M_{2\alpha _{k}}}\left( x\right) =D_{M_{2\alpha _{k}}}\left( x\right)
+\psi _{_{M_{2\alpha _{k}}}}\left( x\right) D_{j}\left( x\right) ,\text{
\qquad when \thinspace \thinspace }j<M_{2\alpha _{k}}.
\end{equation*}

for $II_{2}$ we have:
\begin{eqnarray*}
\left| II_{2}\right| &=&\frac{1}{M}\frac{M_{2\alpha _{k}}}{\alpha _{k}\cdot
q_{\alpha _{k}}}\left| \sum_{j=0}^{q_{\alpha _{k}-1}-1}D_{j+M_{2\alpha
_{k}}}\left( x\right) -D_{M_{_{2\alpha _{k}}}}\left( x\right) \right| \\
&=&\frac{1}{M}\frac{M_{2\alpha _{k}}}{\alpha _{k}.q_{\alpha _{k}}}\left|
\psi _{M_{2\alpha _{k}}}\left( x\right) \sum_{j=0}^{q_{\alpha
_{k}-1}-1}D_{j}\left( x\right) \right| \\
&=&\frac{1}{M}\frac{M_{2\alpha _{k}}}{q_{\alpha _{k}}}\frac{q_{\alpha _{k}-1}%
}{\alpha _{k}}\left| K_{q_{\alpha _{k}}-1}\left( x\right) \right| \\
&\geq &\frac{1}{2M}\frac{q_{\alpha _{k}-1}}{\alpha _{k}}\left| K_{q_{\alpha
_{k}}-1}\left( x\right) \right| .
\end{eqnarray*}

Since

\begin{equation*}
q_{\alpha _{k}}\leq M_{2\alpha _{k}}\left( 1+\frac{1}{4}+...+\frac{1}{4^{n}}%
\right) \leq 2M_{2\alpha _{k}},
\end{equation*}
for $II_{2}$ we obtain

\begin{equation*}
\left| II_{2}\right| \geq \frac{1}{2M}\frac{q_{\alpha _{k}-1}}{\alpha _{k}}%
\left| K_{q_{\alpha _{k}}-1}\left( x\right) \right| .
\end{equation*}

Let $M_{2\alpha _{k-1}+1}-1\leq $ $j<M_{2\alpha _{k}}.$ Then from (\ref{8})
we have

\begin{eqnarray*}
\left| S_{j}f\left( x\right) \right| &=&\left| \sum_{v=0}^{j-1}\widehat{f}%
(v)\psi _{v}\left( x\right) \right| \\
&=&\left| \sum_{v=0}^{M_{2\alpha _{k-1}+1}-1}\widehat{f}(v)\psi _{v}\left(
x\right) \right| \\
&=&\left| \sum_{\eta =0}^{k-1}\sum_{v=M_{2\alpha _{\eta }}}^{M_{2\alpha
_{\eta }+1}-1}\frac{M_{2\alpha _{\eta }}}{M\text{ }\cdot \alpha _{\eta }}%
\psi _{v}\left( x\right) \right| \\
&=&\left| \sum_{\eta =0}^{k-1}\frac{M_{2\alpha _{\eta }}}{M\text{ }\cdot
\alpha _{\eta }}\left( D_{M_{_{2\alpha _{\eta }+1}}}\left( x\right)
-D_{M_{_{2\alpha _{\eta }}}}\left( x\right) \right) \right| \\
&\leq &\frac{2\left( M_{2\alpha _{k-1}}\right) ^{2}}{\alpha _{k-1}}.
\end{eqnarray*}

Hence

\begin{eqnarray}
\left| I\right| &\leq &\frac{1}{q_{\alpha _{k}}}\sum_{j=0}^{M_{2\alpha
_{k}}-1}\left| S_{i}f\left( x\right) \right|  \label{10} \\
&\leq &\frac{2M_{2\alpha _{k}}}{q_{\alpha _{k}}}\frac{\left( M_{2\alpha
_{k-1}}\right) ^{2}}{\alpha _{k-1}}  \notag \\
&\leq &\frac{2\left( M_{2\alpha _{k-1}}\right) ^{2}}{\alpha _{k-1}}.  \notag
\end{eqnarray}

Applying (\ref{4}) we have

\begin{equation*}
\left| I\right| ,\left| II_{1}\right| \leq \frac{2\left( M_{2\alpha
_{k-1}}\right) ^{2}}{\alpha _{k-1}}\leq \frac{1}{16M}\frac{M_{\alpha _{k}}}{%
\alpha _{k}}.
\end{equation*}

Consequently,

\begin{eqnarray}
\left| \sigma _{q_{\alpha _{k}}}f\left( x\right) \right| &\geq &\left|
II_{2}\right| -\left( \left| I\right| +\left| II_{1}\right| \right)
\label{11} \\
&\geq &\frac{1}{8M\cdot \alpha _{k}}\left( 4q_{\alpha _{k}-1}\left|
K_{q_{\alpha _{k}-1}}\left( x\right) \right| -M_{\alpha _{k}}\right) .
\notag
\end{eqnarray}

Denote

\begin{equation*}
I_{2\alpha _{k}}\left( 0,...,x_{2\eta }\neq 0,0,...,0,x_{2s}\neq
0,x_{2s+1},...x_{2\alpha _{k}-1}\right) =I_{2\alpha _{k}}^{\eta ,s}.
\end{equation*}
Let

\begin{equation*}
x\in I_{2\alpha _{k}}^{\eta ,s},\,\,\,\,\,\eta =\left[ \frac{\alpha _{k}}{2}%
\right] ,\left[ \frac{\alpha _{k}}{2}\right] +1,...,\alpha _{k}-3,\,s=\eta
+2,\eta +3,\alpha _{k}-1.
\end{equation*}

Applying lemma 2 we have:

\begin{equation*}
4q_{\alpha _{k}-1}\left| K_{q_{\alpha _{k}-1}}\left( x\right) \right| \geq
M_{2\eta }M_{2s}.
\end{equation*}

Since

\begin{equation*}
2s\geq 2\left[ \frac{\alpha _{k}}{2}\right] +4>\alpha _{k}+1,
\end{equation*}
we have

\begin{equation*}
M_{2s}>M_{\alpha _{k}+1}\geq m_{\alpha _{k}}M_{\alpha _{k}}\geq 2M_{\alpha
_{k}}.
\end{equation*}
Hence

\begin{equation}
M_{2s}M_{2\eta }-M_{\alpha _{k}}\geq \frac{1}{M}M_{2s}M_{2\eta }.  \label{12}
\end{equation}
From (\ref{11}-\ref{12}) we have

\begin{equation*}
\left| \sigma _{q_{\alpha _{k}}}f\left( x\right) \right| \geq \frac{1}{%
8M^{2}\cdot \alpha _{k}}M_{2s}M_{2\eta },\,\,x\in I_{2\alpha _{k}}^{\eta ,s},
\end{equation*}
where

\begin{equation*}
\eta =\left[ \frac{\alpha _{k}}{2}\right] ,\left[ \frac{\alpha _{k}}{2}%
\right] +1,...,\alpha _{k}-3,\,\,s=\eta +2,\eta +3,\alpha _{k}-1.
\end{equation*}

Hence we can write

\begin{eqnarray*}
&&\int_{G_{m}}\left| \sigma _{q_{\alpha _{k}}}f\left( x\right) \right| ^{%
\frac{1}{2}}d\mu \left( x\right) \\
&\geq &\sum_{\eta =\left[ \alpha _{k}/2\right] }^{\alpha _{k}-3}\sum_{s=\eta
+2}^{\alpha _{k}-1}\sum_{x_{2s+1}=0}^{m_{2s+1}-1}...\sum_{x_{_{2\alpha
_{k}-1}}=0}^{m_{2\alpha _{k}-1}-1}\int_{I_{2\alpha _{k}}^{\eta ,s}}\left|
\sigma _{q_{\alpha _{k}}}f\left( x\right) \right| ^{\frac{1}{2}}d\mu \left(
x\right) \\
&\geq &\frac{1}{\sqrt{8}M\sqrt{\alpha _{k}}}\sum_{\eta =\left[ \alpha
_{k}/2\right] }^{\alpha _{k}-3}\sum_{s=\eta +2}^{\alpha _{k}-1}\frac{%
m_{2s+1...}m_{2\alpha _{k}-1}}{M_{2\alpha _{k}}}\sqrt{M_{2s}M_{2\eta }} \\
&\geq &\frac{1}{8M\sqrt{\alpha _{k}}}\sum_{\eta =\left[ \alpha _{k}/2\right]
}^{\alpha _{k}-3}\sum_{s=\eta +2}^{\alpha _{k}-1}\frac{\sqrt{M_{2s}M_{2\eta }%
}}{M_{2s+1}} \\
&\geq &\frac{1}{8M^{2}\sqrt{\alpha _{k}}}\sum_{\eta =\left[ \alpha
_{k}/2\right] }^{\alpha _{k}-3}\sum_{s=\eta +2}^{\alpha _{k}-1}\sqrt{\frac{%
M_{2\eta }}{M_{2s}}}.
\end{eqnarray*}

It is easy to show that

\begin{equation*}
\sum_{s=\eta +2}^{\alpha _{k}-1}\sqrt{\frac{M_{2\eta }}{M_{2s}}}\geq \sqrt{%
\frac{M_{2\eta }}{M_{2\eta +4}}}\geq \frac{1}{M^{2}}.
\end{equation*}

Consequently,

\begin{eqnarray*}
&&\int_{G}\left| \sigma _{q_{\alpha _{k}}}f\left( x\right) \right| ^{\frac{1%
}{2}}d\mu \left( x\right) \\
&\geq &\frac{1}{8M^{2}\sqrt{\alpha _{k}}}\sum_{\eta =\left[ \alpha
_{k}/2\right] }^{\alpha _{k}-3}\left( \sum_{s=\eta +2}^{\alpha _{k}-1}\sqrt{%
\frac{M_{2\eta }}{M_{2s}}}\right) \\
&\geq &\frac{1}{8M^{4}\sqrt{\alpha _{k}}}\sum_{\eta =\left[ \alpha
_{k}/2\right] }^{\alpha _{k}-3}1 \\
&\geq &c\sqrt{\alpha _{k}}\rightarrow \infty ,\text{\qquad as }%
k\longrightarrow \infty .
\end{eqnarray*}

Theorem 1 is proved.

\end{document}